\title{ Singular perturbations and first order PDE on manifolds }
\author {David Holcman \thanks{ Scuola Normale Superiore di Pisa, 7 Piazza dei Cavalieri \, Italy}  \thanks{Weizmann Institute of Science, Rehovot 76100, Israel} and  Ivan Kupka \thanks{ Universit\'e Paris VI , Tour 46, 5 etage,  4 Place Jussieu 75005 Par
is, France} }
\date{ October 1999 }
\documentclass[10pt]{article}

\newtheorem{theorem}{Theorem}
\newtheorem{lem}{Lemma}

\newtheorem{prop}{Proposition}
\font\bb=msbm10 at 12pt
\font\bbbis=msbm10 at 10pt

\def\rR{\hbox{\bb R}} 
\def\rRp{\hbox{\bbbis R}}

\begin{document}
\maketitle
\begin{abstract}
In this note we present some results concerning  the concentration of
sequences of  first eigenfunctions on the limit sets of a Morse-Smale dynamical system on a compact Riemanniann manifold. More precisely a renormalized sequence of eigenfunctions converges to a measure $\mu$ concentrated on the
hyperbolic sets of the field. The set of all possible measure turns out to be a sum of a finite Dirac distributions localized at the critical point of the field  and absolutely continuous measure with respect to the Lebesgue measure on each limit cycles : the coefficients which appear in the limit measure can be characterized using the concentration theory.

In the second part, certain aspects of some first order PDE on manifolds are studied. We study the limit of a sequence solutions of a second order PDE, when a parameter of viscosity tends to zero. Under some explicit assumptions on some vector fields,  bounded and differentiable solutions are obtained. We exibit the role played by the limit sets of the dynamical systems and provide in some cases an explicit representation formula. 
\end{abstract}

\section{Introduction}

$(V_{n},g)$ denotes a compact riemannian manifold of dimension $n \ge
2$, with no boundary and  $\Delta_{g} =-\nabla_{i}\nabla^{i}$ is the
Laplace-Beltrami operator. This note concerns the study of the
operator $ L_{\epsilon} = \epsilon \Delta  +\sum_{ i=1 }^{n}
b_{i}\partial_{i}  + c $ acting on smooth functions, when the
parameter $\epsilon$ converges to zero. $b$ is a regular vector field
and c is a positive function.

In the first part, we are interested to study the behavior of the
first eigenfunction sequence  $u_{\epsilon} $ associated to the the smallest eigenvalues $\lambda_{\epsilon}>0 $ of the operator $ L_{\epsilon}(u_{\epsilon}) = \lambda_{\epsilon}u_{\epsilon} $, when the  parameter $\epsilon$ converges to zero.

Some local results about these sequences are known on bounded domain
of $\rR^{n}$ (see Friedmann \cite{Friedman1} \cite{Friedman2}
\cite{Friedman3}, Friedlin Ventcell' \cite{VentcellF}) : when b has only one attracting point, $u_{\epsilon}$ converges uniformly on every
compact set to a constant as $\epsilon$ converges to 0, but when the point is repulsive the sequence converges in the distribution sense to a Dirac distribution centered at this point.

The behavior of the first eigenvalue $\lambda_{\epsilon}$ as $\epsilon$
goes to zero is well known for a large class of dynamical system : $\lambda_{\epsilon} $ converges to some quantity called the topological pressure P. This number is characterized by a variational problem on the set of probability measures : when the field b has a finite number of hyperbolic invariant sets K, then the topological pressure is 
\begin{eqnarray}
P = \sup_{ } (h_{\mu} +\int_{V_{n}} (c- \frac{d \det D\phi^{u}_{t}}{dt}  ) d\mu \mid
\hbox{ supp } \mu \subset K \,\hbox{ and } \mu \, \phi_{t}-\hbox{invariant} )
\end{eqnarray}
The set of measures considered here is all the measures with support in K
and invariant by the flow. $ h_{\mu}$ is metric entropy (see
\cite{Kifer1}), $\phi_{t}$ is the flow induced by the vector field $b$
and $D\phi^{u}_{t}$ is the differential of $\phi$ restricted to the
unstable bundle of $K$.

The topological pressure P is attained at a measure, called the
equilibrium state (see \cite{Kifer1} \cite{Varadh}). When the recurrence set of the
field $b$ consists of p hyperbolic sets $K_{i}$, $i=1..p$ the measure
$\mu$ is concentrated on the union of the $K_{i}'s$. More precisely, $\mu =  \sum_{i=1}^{p} p_{i} \mu_{K_{i}} $ where  $\mu_{K_{i}}$ is the equilibrium state associate to $K_{i}$, $ p_{i}\geq 0 $ and  $ \sum_{i=1}^{p} p_{i}  =1$  (see theorem 3.4 in \cite{Kifer2}).

Unfortunately this fruitful approach is not adapted to the study
of the first eigenfunction problem since the equilibrium measures are
not invariant by the flow, in general. However we shall  see that the
limit measures have their supports on the recurrent sets of the field
$b$ although they are not invariant by the flow.

In the second part, we gives some results about the behavior of the
solutions of the equation $L_{\epsilon} u_{\epsilon} = f$ when
$\epsilon$ goes to zero, for a  given smooth positive function f. Our main
interest is to find bounded $C^{1}$-solutions and  to understand
the interaction between the geometry of the characteristic curves (trajectories) of the field $b$
and the behaviour of the solutions of $L_{\epsilon} u_{\epsilon} =f$
when $\epsilon$ goes to zero. The fascinating point here is that when
$\epsilon$ goes to zero the elliptic equation   $L_{\epsilon}
u_{\epsilon} = f$ tends to a hyperbolic one. In this last case the
singularities of the solutions propagate along the characteristic
whereas, in the elliptic case there is no propagation of singularities.

The fields  considered here are Morse-Smale ( see \cite{Robinson} \cite{Smale} ), that is 1-) the recurent
set consists of a finite number of hyperbolic stationnary points and
periodic orbits 2-) the stable and unstable manifolds of the recurrent
orbits are pairwise  transversal (for all p,q, the unstable manifold $ W
^{u}(p)$ of p and the stable manifold $ W^{s}(q)$ of q intersect
transversally : $T_{a}( W^{u}(p)) \oplus T_{a}(W^{s}(q)) =T_{a}
(V_{n}) $). A more general class of vector fields will be considered elsewhere. 

\begin{theorem}
Suppose that the first eigenvalue of the operator $\Delta_{g}  +  a  $
is positive  and $a$ is a positive function with a finite number of minimum points which are not degenerate (in the sense of a Morse). Consider the first eigenvalue problem (which has the following variational formulation) 

\begin{eqnarray}
\lambda_{\epsilon} = \inf_{ u \in H_{1}(V_{n})-\{ 0\} } \frac{\epsilon \int_{V_{n}} { \mid \nabla u \mid}^{2}+a u^{2}   } { \int_{V_{n}} u^{2} } 
\end{eqnarray}
Then, when $\epsilon$ converges to zero the sequence
$\lambda_{\epsilon}$ converges to the minimum of the function $a$ and
the set of limits for the weak topology, when $\epsilon$ goes to zero, of the
family of measures $\frac{u^{2}_{\epsilon} dV_{g} }{ \int_{V_{n}}u^{2}_{\epsilon} dV_{g} }$ defined by the positive solution$ u_{\epsilon} $ of the PDE, 
\begin{eqnarray} \label{edpfdt}
\epsilon \Delta_{g} u_{\epsilon} +  a u_{\epsilon}  =  \lambda_{\epsilon} u_{\epsilon}   \mbox{ on } V_{n}\end{eqnarray}
is contained in the simplex $M = \{ \nu =\sum_{i=1}^{n} c_{i}^{2}\delta_{P_{i}}, \sum_{i=1}^{n} c_{i}^{2}= 1 \}$ of all probability measures with support in the finite set $\{ P_{i} \mid i=1..n\}$  where $\delta_{P}$ denotes the Dirac measure at the point $P$.
\end{theorem}
Remark : $u_{\epsilon}$ is uniquely defined up to a multiplicative constant by the Krein-Rutman theorem.

Consider now the case when $b = \nabla \phi $ and the function c is
choosen  so that the eigenvalue $\lambda_{\epsilon}$ of the operator
$ L_{\epsilon} = \epsilon \Delta  +\sum_{ i=1 }^{n} b_{i}\partial_{i}
+ c$ is positive on the manifold. To study the family $u_{\epsilon}$
of solutions of the PDE
\begin{eqnarray} \label{ref5}
\epsilon \Delta u_{\epsilon}  +\sum_{ i=1 }^{n}
b_{i}\partial_{i}u_{\epsilon} + c u_{\epsilon} =\lambda_{\epsilon}
u_{\epsilon} \hbox{ on } V_{n},
\end{eqnarray}
we apply the transformation $ b= \nabla \phi = -2 \epsilon \nabla
\psi_{\epsilon}$ and consider the new variable $ v_{\epsilon} =
u_{\epsilon} \psi_{\epsilon}$. Equation (\ref{ref5}) is transformed
into the following PDE where the vector field disappear:
\begin{eqnarray} 
\epsilon^{2} \Delta v_{\epsilon}+ a_{\epsilon} v_{\epsilon} = \epsilon
\lambda_{\epsilon} v_{\epsilon} \hbox{ on } V_{n}
\end{eqnarray}
where $a_{\epsilon} = c\epsilon +\frac{ { \mid \nabla \phi  \mid}^{2}}{4}
+\frac{\epsilon \Delta \phi}{2}.$
\begin{prop}
Suppose that the following condition is satisfied : at each minimum points $P_{i}$ of the function $c$, $c(P_{i})+ \Delta \phi(P_{i})/2 \geq 0$.
Let $v_{\epsilon}$ be a minimizer of the following variational problem
\begin{eqnarray} 
\epsilon \lambda_{\epsilon} = \inf_{ u \in H_{1} (V_{n})-\{ 0\} } \frac{\epsilon^{2} \int_{V_{n}} { \mid \nabla u  \mid}^{2}+ a_{\epsilon} u^{2}   } { \int_{V_{n}} u^{2} } 
\end{eqnarray}
then 
\begin{itemize}
\item $\lim_{\epsilon \rightarrow 0} \lambda_{\epsilon} = \inf_{V_{n}}  (\nabla \phi)^{2} =0$
\item The set of limits for the weak topology of the measures  $v^{2}_{\epsilon}dV_{g}$ is localized at the points where the function $
  \lim_{\epsilon \rightarrow 0} \min  a_{\epsilon}  $ is
  zero. These measures have support on the critical
  points of the vector field $b$.

 Moreover $ \sup_{ V_{n} } v_{\epsilon}$ converges to $ + \infty$
\end{itemize}
\end{prop}
Going back to the sequence of eigenfunctions, we get that for some subsequence of $\epsilon$ tending to zero
\begin{eqnarray} 
\lim_{\epsilon \rightarrow 0 }\frac{ \int_{V_{n}} e^{-\phi/\epsilon}
  u^{2}_{\epsilon} a_{\epsilon}  }{\int_{V_{n}} e^{-\phi/\epsilon}
  u^{2}_{\epsilon}  } = \min (\nabla \phi)^{2} =0
\end{eqnarray}
and for any function  $\psi \in C^{\infty}(V_{n})$ we have for some subsequence of  $\epsilon$ tending to zero
\begin{eqnarray} 
\frac{ \int_{V_{n}} e^{-\phi/\epsilon} u^{2}_{\epsilon} \psi  }{\int_{V_{n}} e^{-\phi/\epsilon} u^{2}_{\epsilon}  } \rightarrow \sum_{i=1}^{n} c^{2}_{i}\psi(P_{i})
\end{eqnarray}
where $ c_{i}^{2} = \lim_{\epsilon \rightarrow 0} \frac{ \int_{B_{P_{i}}(\delta) }  e^{-\phi/\epsilon} u^{2}_{\epsilon}}{ \int_{V_{n} }  e^{-\phi/\epsilon} u^{2}_{\epsilon}} $ (the limit is independant of $\delta$).
The set of weak limits of the family of measure  $\frac{
  e^{-\phi/\epsilon} u^{2}_{\epsilon}}{ \int_{V_{n} }
  e^{-\phi/\epsilon} u^{2}_{\epsilon}} dV_{g}$ is given by
$\sum_{i=1}^{n} c^{2}_{i}\delta_{P_{i}}$ where $P_{i}$ are  the critical points of the function  $\phi$, the zero of the vector field $b=\nabla \phi $.

Next we will consider  vector fields of the following form $ b = -
\nabla L + \Omega $ where $ \Omega$ is fixed but not of  gradient type
and we will construct  $L $ to be a type of  Lyapunov function of $
\Omega$. Using this vector field, we study the behavior  of the
sequence of the first eigenfunction for the operator $L_{\epsilon}$
and we obtain the following  theorem  :

\begin{theorem}
On a compact Riemannian manifold $V_{n}$, consider a Morse-Smale vector field $ \Omega$ not a gradient whose recurrent set consists of the stationnary points $P_{1},...,P_{M}$ and of the periodic orbits $\Gamma_{1}$,.., $\Gamma_{N}$. $L$ denotes a special Lyapunov function associated with $\Omega$ defined in the proof of Theorem 2.

For $\epsilon >0 $ let $\lambda_{\epsilon}$, $ u_{\epsilon} $ denotes respectively the first eigenvalue and an associated eigenfunction of the operator 
\begin{eqnarray}
\epsilon \Delta_{g} + b \nabla  +  c   \mbox{ on } V_{n} 
\end{eqnarray}

Then the  set of all weak limits as $\epsilon$ goes to zero, of the family of normalized measures  
\begin{eqnarray}
\ \frac{  e^{-L/\epsilon} u_{\epsilon}^{2}  dV_{g} } {\int_{V_{n}} e^{-L/\epsilon} u_{\epsilon}^{2} dV_{g} }
\end{eqnarray}
is contained in the set $ M = \{ \nu \mid \nu  \hbox{ Borel  measure, support $\nu$ }\subset \cup_{j=1}^{N} \Gamma_{j} \cup \{ P_{i} \mid 1\leq i \leq M \} \hbox{ and } \, \nu = \sum_{j=
1}^{N} f^{2}_{k}(l_{j})dl_{j} +  \sum_{i=1}^{M}  c_{i}^{2} \delta_{P_{i}} \} $ where $dl_{j}$ denotes the arc length on $\Gamma_{j}$, for the metric $g$, $\delta_{P_{i}}$ the Dirac measure with support $P_{i}$, $f_{j}: \Gamma_{j}\rightarrow \rR$ are continuous functions $ 1 \leq j \leq N$ and the $c_{i}$ are constants.

\end{theorem}

To prove the theorem, we will need 3 lemma

\begin{lem} 
There exists a Lyapunov function L for the field $\Omega$, such that
 $L$ is twice differentiable in the neighborhood of the recurrent set of $\Omega$
 and reaches its minimum  on the union of the neighborhoods on the
 recurrent sets only. Outside these neighborhoods, L can
 be arbitrary such that $\Psi(L) = \frac{1}{4}({\mid \nabla L \mid}^{2} +2 (\nabla L,\Omega)) \geq 0  $
\end{lem} 
The same type of Lyapunov function $L$ was constructed by Kamin
\cite{Kamin} \cite{Kamin2} in the case of an attractive point.

\begin{lem}
Under the assumptions of Theorem 2 on the vector field $\Omega$, 
\begin{eqnarray}
\lim_{\epsilon \rightarrow 0} \epsilon  \lambda_{\epsilon} = 0 = \min_{V_{n}} \Psi
\end{eqnarray}
where $\Psi = \frac{(\nabla L)^{2}}{4} +\frac {(\Omega, \nabla L)}{2} $.
\end{lem}

\begin{lem}
Under the assumptions of Theorem 2  on the vector field $\Omega$,
all weak limits of measures $v_{\epsilon}^{2}dV_{g}$ as $\epsilon$ goes to zero are concentrated on the minimum set of the function $\Psi$.
\end{lem}

The measure $ \mu $ is absolutely continuous with respect to
the measure induced by the length along the periodic orbit. This is
true  on each pairwise orbit of the vector field. We have

\begin{eqnarray}
\frac{ d\mu}{dl}= f^{2}(l)= \lim_{\epsilon \rightarrow 0} \int_{ H_{l} } u^{2}_{\epsilon} d\Sigma_{g}
\end{eqnarray}

$H_{l}$ denotes any hypersurface cutting the orbit at the point of abscisse $l$ transversally.
Remark: The limit is independant of the choice of the hypersurface.

\section{ First order PDE on manifolds  }
We study the limit of the solution  $u_{\epsilon}$  of $L_{\epsilon}
(u_{\epsilon}) =f$, where $f$ is a given smooth function on the
manifold, when $\epsilon$ tends to zero. The limit of the sequence
when $\epsilon$ goes to zero,  solves some first order PDE. For some
previous works see \cite{Sacker}, \cite{Sacker2} and \cite{Moser}.
$ c$ and $f$ will denote  two given positive smooth functions and b a
vector field. 
We suppose that $ c_{0}= \inf_{V_{n} } c >0 $ and $b_{0} = \sup_{X \in
  TV_{n} } 1/2(\nabla_{i} b^{k}+\nabla_{k} b^{i})(X_{i},X_{k}) $ a
finite number. S denotes the set of separatrices associates to the
dynamical system.  We prove the following theorem :

\begin{theorem} \label{theex}
On a compact Riemannian manifold, consider a Morse-Smale vector field
$b$ and let $c$ be a positive function satisfying $c(x)\geq c_{0}>0$ and $ c_{0} - b_{0}> 0 $. f is a differentiable function.
Under these assumptions,  there exists a solution $ u \in C^{0}(V_{n})$ such that $ \mid \nabla u \mid \in L_{\infty}(V_{n})$ of the first order PDE :
\begin{eqnarray}
 <b, \nabla u>  +  c u  = f \hbox{ on } V_{n} 
\end{eqnarray}
Moreover if the recurrent set of the vector field $b(x)$ consist of a
finite number of points $P_{1},..P_{p}$, and $c=c_{0}$ is constant
larger than the eigenvalues of $Db$, then u is in $C^{1}(V_{n}-S)$ and
$u$ is unique and is completely determined by the values   
$u(P_{i})=\frac{ f(P_{i})}{ c(P_{i}) }$.
When $b$ possess some limit cycles, the solution is not unique and has no limit near the limit cycles.
\end{theorem}
We adapt the previous theorem to the nonlinear first order PDE on a compact manifold.
\begin{eqnarray} \label{edpnl}
<b(u,x),\nabla u >  + c(u,x)u  = f
\end{eqnarray}
where $b(\lambda,x)$ is a regular vector field, $\lambda$ is parameter  and  $c(\lambda,x)$ and $f$ are two given functions.
The purpose of this part is to find  some conditions on f,c and b to insure
regular solutions, since there exists examples where shocks occur. 
To prove the existence of solutions for (\ref{edpnl}),  we use an
elliptic regularization  and proceed by successive approximations, (see Jausselin
et al. \cite{Moser2}).

For the last theorem we need the following notations :
\begin{eqnarray} 
b_{0} =\frac{1}{2} \sup_{\parallel X \parallel =1, x \in V_{n}, \lambda \in \rRp }( \nabla_{i} b_{k}(x,\lambda) +\nabla_{k} b_{i}(x,\lambda )) (X_{i},X_{k}) \nonumber
\end{eqnarray} 
and $ \gamma = \sup_{\lambda \in \rRp,x \in V_{n}} \mid\partial_{\lambda} b(\lambda,x)\mid $, 
$a_{0} = \inf_{\lambda,x} c(\lambda,x) - b_{0}$.
$A =  \sup_{V_{n}} \mid \nabla f \mid +\sup_{V_{n}} (f/c) \times
  \sup_{V_{n}\times \rRp } \mid \frac{\partial c }{\partial x} \mid $
  and $ \beta = \sup_{V_{n}} \mid c'\mid. \sup_{V_{n}} \frac{f}{c}$, 
  where we use the notation $c' =\frac{ \partial c}{\partial \lambda}$.

We make the following assumptions for the rest of this section
\begin{enumerate}
\item $a_{0}> \beta $ that is  $ \inf_{\lambda,x} c(\lambda,x) - b_{0} >  \sup_{V_{n}} \mid c'\mid \sup_{V_{n}} \frac{f}{c}   $
\item $ (\inf_{\lambda,x} c(\lambda,x)  - b_{0})^{2} +(\sup_{V_{n}} \mid
c'\mid \sup_{V_{n}} \frac{f}{c})^{2} \geq 2(\inf_{\lambda,x}
c(\lambda,x)  - b_{0})(\sup_{V_{n}} \mid c'\mid) \sup_{V_{n}} \frac{f}{c} + $ 
$ 4 \sup_{\lambda,x} \mid \partial_{\lambda} b(\lambda,x) \mid \times \sup_{V_{n}} \mid
\nabla f \mid $

\item If $\Lambda = \inf_{ \rRp \times \rRp } (c(\lambda,u) +uc'(\lambda,u) - b_{0})$ and  $ \Lambda^{2} - 4 A \gamma \geq 0$
\end{enumerate}
These assumptions mean that the minimum of the function c must be
large enough with respect to the vector field b : it is a hyperbolicity condition.
\begin{theorem} \label{theexp}
On a compact manifold, consider the parametrized smooth vector field $b(\lambda,x)$, and  let $c(\lambda,x)$ satisfie $c(\lambda,x)\geq c_{0}>0$ and $c_{0}$ large so that  $c_{0}( c_{0}-b_{0}) >\sup_{V_{n}} (f\mid c' \mid) $ and the conditions above satisfied.
Then there exists a solution $ u \in W^{\infty,1}(V_{n})$ such that the first order PDE :
\begin{eqnarray}
 <b(u,x), \nabla u>  + c(u,x)  u  = f \hbox{ on } V_{n} 
\end{eqnarray}
Moreover if the limit set of the vector field $b(u,x)$ is  a union of
a finite  number of points $P_{1},..P_{p}$, $u$ is in $C^{1}(V_{n}-S)$, where $S$ is the set of separatrices of the  field $b(u,x)$ .  If the number of solution of $ c(u(P_{i}),P_{i} )u(P_{i
})=f(P_{i})$ is finite and equal to k, the number of solutions $u$ is $k^{p}$.
\end{theorem}


\end{document}